\documentclass[a4paper,10pt]{article}
\usepackage{amsmath}
\usepackage{amsfonts}
\usepackage{amssymb}
\usepackage{array,longtable}
\usepackage{amscd}
\usepackage[cp1251]{inputenc}
\usepackage[english]{babel}
\usepackage[usenames]{color}
\usepackage{colortbl}

\begin{document}
\newtheorem{theorem}{Theorem}

\begin{center}
{\large {\bf On matrices with $Q^2$-scalings}} \\[0.5cm]
\end{center}

\begin{center}
    {\bf  O.Y. Kushel}
\end{center}
\begin{center}
  {\it  kushel@mail.ru}
\end{center}
\begin{center}
    Institut f\"{u}r Mathematik, MA 4-5, Technische Universit\"{a}t Berlin, \\ D-10623 Berlin, Germany
\end{center}

\begin{abstract}
We provide a counterexample to some statements dealing with a sufficient property for the square of a matrix to be a $P_0^+$-matrix.

$P$-matrices; $Q$-matrices; $P$-matrix powers
 
Primary 15A48  Secondary 15A18  15A75
 \end{abstract}

Let us recall the following definitions and notations (see, for example, \cite{HERK2}, \cite{PINK}). If $\mathbf A$ is an $n \times n$ matrix, ${\mathbf A}^{(j)}$ ($1 \leq j \leq n$) denotes its $j$th compound matrix, i.e. the matrix which consists of all the minors of the $j$th order of $\mathbf A$, numerated lexicographically.

An $n \times n$ matrix $\mathbf A$ is called a {\it $Q$-matrix} if its sums of principal minors of the same order are all positive (this is equivalent to the following conditions: ${\rm Tr}({\mathbf A}^{(j)}) > 0$ for all $j = 1, \ \ldots, \ n$). An $n \times n$ matrix $\mathbf A$ is called a $P_0$- ($P_0^+$-) matrix if all its principal minors are nonnegative (respectively, nonnegative with at least one positive principal minor of each order). An $n\times n$ matrix $\mathbf A$ is called {\it anti-sign symmetric} if it satisfies the following conditions:
$$A \begin{pmatrix} \alpha \\ \beta \end{pmatrix}A \begin{pmatrix} \beta \\ \alpha \end{pmatrix}  \leq 0 $$
 for all sets of indices $\alpha, \beta \subset \{1, \ \ldots, \ n\},$ $\alpha \neq \beta$, $|\alpha| = |\beta|$.

The following statement was claimed to be proven in \cite{HERK2} (see \cite{HERK2}, p. 115, Proposition 4.4).
\begin{theorem}
Let $\mathbf A$ be a square matrix. If for every positive diagonal matrix $\mathbf D$ the matrix $({\mathbf D}{\mathbf A})^2$ is a $Q$-matrix then ${\mathbf A}^2$ is a $P_0^+$-matrix.
\end{theorem}

This statement does not hold. Let us consider the following counterexample.

{\bf Counterexample.} Let $${\mathbf A} = \begin{pmatrix}
1 & 2 \\ -1 & 5
 \end{pmatrix}. \eqno(1)$$
In this case, we have
$${\mathbf A}^{(2)} = \det({\mathbf A}) = 7.$$

Multiplying by an arbitrary positive diagonal matrix ${\mathbf D} = {\rm diag}\{d_{11}, \ d_{22}\}$, we obtain:
$$ {\mathbf D}{\mathbf A} = \begin{pmatrix} d_{11} & 2d_{11} \\ -d_{22} & 5d_{22} \end{pmatrix};$$
$$ ({\mathbf D}{\mathbf A})^2 = \begin{pmatrix} d_{11}^2 - 2d_{11}d_{22}  & 2d_{11}^2 + 10 d_{11}d_{22} \\
-d_{11}d_{22} - 5d_{22}^2 & -2d_{11}d_{22} + 25d_{22}^2 \end{pmatrix};$$
$$ (({\mathbf D}{\mathbf A})^2)^{(2)} = \det(({\mathbf D}{\mathbf A})^2) = 49 d_{11}^2 d_{22}^2.$$

It is easy to see that
$$ {\rm Tr}(({\mathbf D}{\mathbf A})^2) = d_{11}^2 - 4d_{11}d_{22} + 25d_{22}^2 = (d_{11} - 2d_{22})^2 + 21d_{22}^2 > 0;$$
$$\det(({\mathbf D}{\mathbf A})^2) = 49 d_{11}^2 d_{22}^2 > 0$$
 for any positive values $d_{11}, \ d_{22}$. Thus the matrix $({\mathbf D}{\mathbf A})^2$ is a $Q$-matrix for every positive diagonal matrix $\mathbf D$.
 However,
 $${\mathbf A}^2 = \begin{pmatrix} -1 &  12 \\
 -6 & 23 \end{pmatrix}$$
 is not even a $P_0$-matrix since it has a negative entry on the principal diagonal.

The flaw in the proof is as follows. For a given proper subset $\alpha$ of $\{1, \ \ldots, \ n\}$, the authors construct a positive diagonal matrix ${\mathbf D}_\epsilon$:
$$({\mathbf D}_\epsilon)_{jj} = \left\{
\begin{array}{cc}
1, & j \in \alpha  \\
\epsilon, & j \not\in \alpha \\
\end{array}\right.$$
and claim the following equality for the principal minors: $({\mathbf D}_0{\mathbf A})^2[\alpha] = {\mathbf A}^2[\alpha]$. However, this is not true. $({\mathbf D}_0{\mathbf A})^2[\alpha]$ gives the determinant of $({\mathbf A}_{\alpha})^2$ where ${\mathbf A}_{\alpha}$ is a principal submatrix of $\mathbf A$ spanned by rows and columns with the numbers from $\alpha$, while ${\mathbf A}^2[\alpha]$ gives the determinant of the corresponding submatrix of ${\mathbf A}^2$ (note, that $({\mathbf A}_{\alpha})^2 \neq ({\mathbf A}^2)_{\alpha}$). For example, if $n =3$, $\alpha = \{1, \ 2\}$, ${\mathbf A} = \{a_{ij}\}_{i,j = 1}^3$, we have ${\mathbf D}_{\epsilon} = {\rm diag}\{1,1,\epsilon\}$ and
$${\mathbf D}_0{\mathbf A} = \begin{pmatrix} a_{11} & a_{12} & a_{13} \\ a_{21} & a_{22} & a_{23} \\ 0 & 0 & 0 \\ \end{pmatrix}.$$ In this case, $({\mathbf D}_0{\mathbf A})^2[1,2] = (a_{11}a_{22}-a_{21}a_{12})^2 = \left( A\begin{pmatrix} 1 & 2 \\ 1 & 2 \end{pmatrix}\right)^2$, which is always positive. However, ${\mathbf A}^2[1,2]$ is equal to $\left( A\begin{pmatrix} 1 & 2 \\ 1 & 2 \end{pmatrix}\right)^2 + A\begin{pmatrix} 1 & 2 \\ 1 & 3 \end{pmatrix}A\begin{pmatrix} 1 & 3 \\ 1 & 2 \end{pmatrix} + A\begin{pmatrix} 1 & 2 \\ 2 & 3 \end{pmatrix} A\begin{pmatrix} 2 & 3 \\ 1 & 2 \end{pmatrix}$ and obviously in general case is not equal to $({\mathbf D}_0{\mathbf A})^2[1,2]$.

The following statements were claimed to be proven in \cite{HERK2} using false Proposition 4.4 (see \cite{HERK2}, p. 115, Proposition 4.6 and p. 116, Theorem 4.8).

\begin{theorem} Let $\mathbf A$ be a $2 \times 2$ matrix. Then the following are equivalent:
\begin{enumerate}
\item[{\rm (i)}] For every positive diagonal matrix $\mathbf D$ the matrix $({\mathbf D}{\mathbf A})^2$ is a $Q$-matrix.
\item[{\rm (ii)}] The matrix ${\mathbf A}^2$ is a $P_0^+$-matrix.
\end{enumerate}\end{theorem}

\begin{theorem} Let $\mathbf A$ be an anti-sign symmetric matrix. Then the following are equivalent:
\begin{enumerate}
\item[{\rm (i)}] For every positive diagonal matrix $\mathbf D$ the matrix $({\mathbf D}{\mathbf A})^2$ is a $Q$-matrix.
\item[{\rm (ii)}] The matrix ${\mathbf A}^2$ is a $P_0^+$-matrix.
\end{enumerate}\end{theorem}
The implication $(i) \Rightarrow (ii)$ is false in both of the statements. An anti-sign symmetric $2 \times 2$ matrix $\mathbf A$ given by Formula (1) provides the counterexample for both of them. Thus we conclude that Proposition 4.4 fails even in the case of $2 \times 2$ matrices.

\end{document}